# An Interactive Multi-Dimensional Flexibility Scheduling in Low-carbon Low-inertia Power Systems

Farhad Samadi Gazijahani, *Member, IEEE*, Rasoul Esmaeilzadeh, *Member, IEEE*

*Abstract*—Today, electrical energy plays a significant and conspicuous role in contemporary economies; as a result, governments should place a high priority on maintaining the supply of electrical energy. In order to assess various topologies and enhance the security of power systems, it may be useful to evaluate robustness, dependability, and resilience all at once. This is particularly true when there is a significant amount of renewable energy present. The R3 concept, which consists of these three interrelated characteristics, describes the likelihood that a power system would fail, the potential severity of the repercussions, and the speed at which the system will recover from a failure. This paper uses eight case studies created from the IEEE 24-bus RTS and thoroughly assesses the properties of reliability, robustness, and resilience to highlight the significance of the issue. The sequential Monte Carlo method is used to evaluate reliability, cascade failure simulations are used to evaluate robustness, and a mixed-integer optimization problem is used to study resilience. Different indicators related to each of the three assessments are computed. The significance of the combined analysis is emphasized as the simulation findings are described visually and statistically in a unique three-dimensional manner eventually.

## I. Introduction

Electrical power systems must be dependable, strong, and resilient to handle both expected and unanticipated circumstances. If the system is adaptable enough, these objectives can be accomplished. The capacity of the system to deploy its resources in order to respond to a disruption (or changes in net load) enough for maintaining the system's security may be referred to as operational flexibility. Generation sources and energy storage systems, which are often contracted to supply regulating energy (spinning reserves) or are able to be re-dispatched quickly enough, are the main sources of flexibility (manual reserves). The flexibility may be quantitatively assessed in terms of available energy, power capacity, and up- and down-ramping, it is important to note.

In addition, electricity is becoming increasingly crucial for the everyday and continuous operation of daily activities in the present decarburization process of contemporary civilizations, particularly in the low-carbon low-inertia power systems. As a consequence, threats and interruptions to the security of the electrical system are growing and changing at a similar pace. More research is thus required to examine the related characteristics, guarantee that systems are continually becoming more secure, and monitor how systems are evolving in response to various threats to the industry.

The scientific literature [1] clearly defines the differences between the notions of dependability, robustness, and resilience in a power system. According to Georges Zissis' message in the IEEE Industry Applications Magazine [2], dependability is the likelihood that a system will operate well for a certain amount of time when it is employed in a specific operating environment. In the case that one or more assets are lost, this feature assesses the network performance. On the other side, robustness refers to the degree to which a system can continue to operate even when one or more of its components or subsystems fail. This feature, which is more assertive than dependability, takes the loss of many assets into account and measures network performance in the case of cascading failures. The last definition of resilience is a system's capacity to recover from a significant disturbance within acceptable degradation boundaries. Generally speaking, this idea examines High Impact Low Probability (HILP) occurrences, such as major natural catastrophes [3]. The "R3 idea" is the current name for the three joint qualities [2]. R3 is shown as an example in Fig. 1.

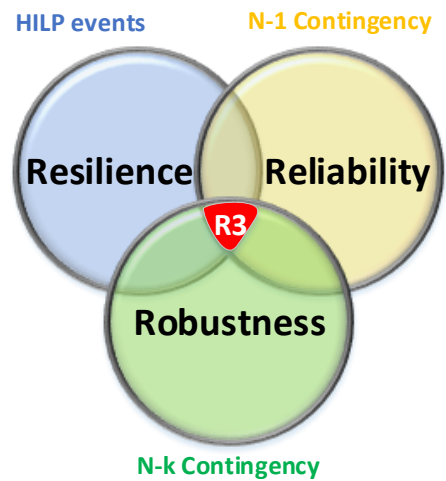

**Fig. 1. Schematic representation of the R3 concept [4].**

Today, there is a significant demand to boost an electrical network's functionality and quality. This desire is a result of society's growth and change toward more robust, sustainable, and carbon-free societies. The reliability, robustness, and resilience (R3) concepts are a body of research that calls for the development of fresh, integrated methodological frameworks to investigate their many linked properties. Some of these traits may be studied using techniques that are described in scholarly literature [5]. For instance, earlier research has used a variety of techniques, including energy hub-based techniques, model order reduction, metaheuristic searching genetic algorithms, multi-criteria decision analysis, advanced intelligent strategies, as well as linear programming [6–11].

Cascading failures, on the other hand, provide one of the largest difficulties in researching the R3 idea. These occurrences are difficult to understand since there are many possible causes or origins, making it almost impossible to address them all. In this regard, complex network theory (also known as graph theory) may be appropriate to describe dynamic behavior, examine disturbance propagation, and determine how structurally resistant a system is [12–15]. It should be mentioned that this method may also be used to locate crucial assets and the circumstances that could lead to a cascade failure [16–21].

Some studies suggest that as weather conditions have the potential to have a significant influence on the system and the sectors connected to it, reliability and robustness studies should take them into consideration. These studies tend to be more focused on extreme disruptions. According to this assumption, various articles provide measurements, defense tactics, and theoretical frameworks to thoroughly investigate this issue [22–24].

Resilience is a brand-new field of study that includes methods and strategies for dealing with issues related to safeguarding and restoring the functions offered by a power system. Some researchers assess a network's resilience by taking into account its topology's strengths and weaknesses as well as its ability to transmit power during various disturbances including earthquakes, floods, and natural catastrophes [25–27]. Note that it is crucial to assess a network's resilience after a high-impact disturbance due to the rise in interruptions [28]. To design the system's iterative recovery after a disruptive event, several researchers have suggested approaches [29–31]. These studies help to increase a power system's responsiveness.

Changes in the generating mix of power systems are taken into account in several combined robustness and resilience studies. As an example, some studies examine cutting-edge techniques while others provide comprehensive analyses of decision-making [32, 33]. The trait of resilience is also covered by other research, which include definitions, measurements, recommendations, practical difficulties, and technological issues [34–37].

The scientific literature has works that explore reliability, robustness, and resilience from a variety of angles, but there aren't many research that take all three ideas into account simultaneously. The development of a theoretical and data-based methodological framework to examine the features and interactions between all concepts in an electrical power system is the particular purpose of this research, which is motivated by the latter. The results of the reliability, robustness, and resilience studies taken together may provide a more accurate picture of network performance. Future study in this field may be encouraged by include the three features in a combined analysis. It is crucial to note that this research instead highlights the need of taking these notions' crucial roles for the electrical network into account in an integrated way rather than discussing how to enhance the study of these features. Decision-makers should prioritize improving energy security since a dependable power system may not be strong or resilient to other threats or disruptions.

Applying the sequential Monte Carlo technique, the reliability assessment is done for all of the aforementioned factors and based on the provisions and guidelines found in the scientific literature. The reliability assessment measures the Expected Energy Not Supplied (EENS), the Expected Demand Not Supplied (EDNS), the Expected Frequency of Load Curtailment (EFLC), the Loss of Load Expectancy (LOLE), the Loss of Load Probability (LOLP), and the Average (ADLC). The robustness assessment is performed by modeling cascade failures and estimating the Satisfied Demand (SD) at each level of system disintegration. This iterative approach eliminates an asset, quantifies the power flows in the network, removes the system's overloaded linkages, then detects and balances the resultant subsystems to decide whether the cascading event continues or terminates. The proposed mixed-integer optimization problem for the resilience research uses real and integer variables to represent the planned dispatch of the generators and the operational condition of the power lines, respectively. This technique identifies the power lines that need to be closed and the re-dispatch of generating plants for the best possible recovery of network connection using the disintegration condition of the system at the conclusion of the cascading failure as input data. At each step of the recovery process, the Recovered Demand (RD) is determined, and the Energy Not Supplied (ENS) is measured at the conclusion.

The three preceding approaches employ DC (direct current) power flows as it offers rapid results. Other techniques may, of course, be used depending on the level of precision necessary for the findings, but because the main goal in this case is to create an integrated framework, this model is adequate. In a sequential research framework, the various outcomes are examined visually and quantitatively. After that, the combined analysis of the three ideas is offered. The suggested method might have a substantial influence on a power network's performance and quality, increase customer happiness, and assist planners in making better investment decisions for network topologies. Eight case studies based on the IEEE 24-bus RTS are subjected to numerical testing to examine the similarities and differences between the ideas [38].

The structure of this article is as follows: The dependability, robustness, and resilience processes used to assess a power system's performance in a structured and methodical way are described in Section 2. Case studies based on the well-known IEEE 24-bus reliability test system are presented in Section 3. The simulation findings that were produced after using the methods outlined above are discussed in Section 4. Section 5 concludes by summarizing the key findings of this study project.

**II. Reliability, robustness and resilience methodologies**

The methods for assessing an electrical power system's dependability, robustness, and resilience are detailed in this section. In general, the sequential Monte Carlo approach is used for the reliability assessment, cascade failure simulations are used for the robustness evaluation, and mixed-integer optimization problems are created for the resilience evaluation. These three techniques adhere to the principles outlined in the academic literature.

### 2.1. Reliability procedure

Adequacy and security are two categories of reliability that are well defined in the scientific literature. Adequacy measures, on the one hand, how well the generating capacity adapts to the demands and limitations of the system. Security, on the other hand, examines a power system's resilience to the failure of one or two components. The safety of an electrical system is the main subject of this study. This kind of assessment may be carried out from either an analytical or a simulation standpoint. The first one needs beginning assumptions to simplify the issue and create an analytical model, which might cause the analysis that results to become irrelevant. The second one replicates the system's random behavior via several trials and takes into account all potential network variables. This study endeavor uses the Monte Carlo technique because it uses a simulation-based methodology [39, 40]. The sequential and organized phases of this process are detailed in Algorithm 1.

---

**Algorithm 1** Reliability

**Input:** technical data of the power system and iterations ($N$).
**Output:** statistical indicators EENS, EDNS, EFLC, LOLE, LOLP and ADLC.
**Step 1.** *Start*: establish the operational state of the assets, that is, normal or failure.
**Step 2.** *Modelling of outages*: these events are modelled using the meantime to failure (MTTF) and the meantime to repair (MTTR). These indicators are inversely related to failure ($\lambda$) and repair ($\mu$) rates of the assets,

$$MTTF = \frac{1}{\lambda}; \ Unit:[h] \quad (1)$$

$$MTTR = \frac{1}{\mu}; \ Unit:[h] \quad (2)$$

**Step 3.** Time between states: quantify the time that the assets spend in normal and failure state, that is, the time to repair (TTR) and the time to failure (TTF),

$$TTR = -\ln(r) \times MTTR; \ Unit:[h] \quad (3)$$

$$TTF = \frac{-\ln(r)}{\lambda} \times 8760; \ Unit:[h] \quad (4)$$

Where $r$ is a random number uniformly distributed between [0,1]. This step is repeated for a specific time, usually one year.
**Step 4.** *Overlapping time*: calculate the overlapping times of failures of the elements (when several components are simultaneously out of service) for a temporal resolution of 1 hour in a time horizon of 1 year, that is, 8760-time steps of 1 hour each.
**Step 5.** *Power flows*: carry out a power flow study considering the operational state of the components throughout the year.
**Step 6.** *Reliability indicators*: evaluate the security of the power system by reliability indices (5)-(10), using the results from Step 5,

$$EENS = \frac{\sum_{i=1}^{N_y} \sum_{j=1}^{N_i} E_{j,i}}{N_y}; \ Unit:[MWh/yr] \quad (5)$$

Where $E_{j,i}$ is the energy not supplied in disruption $j$ of iteration $i$, $N_y$ is the total number of simulated years and $N_i$ is the total number of disruptions in a year $i$.

$$EDNS = \frac{EENS}{8760}; \ Unit:[MW] \quad (6)$$

$$EFLC = \frac{\sum_{i=1}^{N_y} N_i}{N_y}; \ Unit:[outages/yr] \quad (7)$$

$$LOLE = \frac{\sum_{i=1}^{N_y} \sum_{j=1}^{N_i} D_{j,i}}{N_y}; \ Unit:[h/yr] \quad (8)$$

Where $D_{j,i}$ is the duration of disruption $j$ in iteration $i$.

$$LOLP = \frac{LOLE}{8760}; \ Unit:[\%] \quad (9)$$

$$ADLC = \frac{LOLE}{EFLC}; \ Unit:[h/outages] \quad (10)$$

**Step 7.** *Iterations*: repeat the previous steps until obtaining a covariance less than 5% for the EENS index [41].

---

This technique generally presupposes that each electrical system component may exist in one of two states: operational or inoperative. Then, it makes the assumption that the component's residence time has an exponential distribution and that the state transition is influenced by both the current state and the transition rates. The failure and repair rates of the components are the transition rates between the two states. Each component's random samples of its state are statistically dependent, or connected to the preceding sample. After determining the overlapping periods, it next performs power flows and computes the reliability indicators for the electrical system under investigation. This process is performed several times until the EENS indicator's covariance is less than 6%.

### 2.2. Robustness procedure

Blackouts are often the result of a cascade failure, which is a succession of interconnected events such voltage issues and the disconnecting of electrical cables and loads. Modeling these occurrences is challenging since they include a large number of very dynamic events. However, since they have a significant impact on millions of individuals and result in significant economic losses, it is crucial to investigate and model them [12]. This paper measures resilience both before and after a cascade failure in operational areas [42]. Here, the performance of an electrical power system under such disruptions is gauged using the SD index. This index, which ranges from 1 to 0, measures the assets that are isolated during a network's breakdown. The effect on unconnected loads grows as the SD index lowers. The sequential and methodical procedures used in Algorithm 2 to mimic cascading failures in an electrical power system are described.

---

**Algorithm 2** Robustness

**Input:** technical data of the power system and overload tolerance parameter of the lines ($\alpha$).
**Output:** degradation of the electrical power system. Satisfied Demand (SD) in each disintegration stage ($s$), set of subnets ($I$), set of isolated elements ($E$) and state of the links ($\mu_k$), i.e. open or closed.

**Step 1.** *Start*: $SD_{base}=D_{load}$, $I=\{0\}$ and $E=\{0\}$. At the beginning, all the lines are operational.

**Step 2.** *DC power flows*: calculate the flows ($P$) in each link ($k$) within the network and determine the maximum power threshold ($P_k^{max}$) of the lines using (11),

$$P_k^{\max} = \alpha_k \times P_k ; \; Unit : [MW] \qquad (11)$$

**Step 3.** *Initiating event*: randomly remove an asset from the system. The latter represents the event that triggers the cascading failure.

**Step 4.** *Increase or decrease flows*: determine the increases or decreases in each power line; initialize $s=1$ as the first disintegration stage.

**Step 5.** *Triggering of switches*: evaluate the condition $|P_k^s| < P_k^{max}$ in all power lines of the system. Remove all overloaded links, i.e. $|P_k^s| > P_k^{max}$, and go to Step 6; otherwise, go to Step 10.

**Step 6.** *Transversal graph algorithm*: use the DFS algorithm to determine the subnets ($I$) and isolated elements ($E$) formed after the activation of the switches.

**Step 7.** *Energy balance:*
  a) for each island $I_i$ with generators $g \in I_i$ evaluate,
  - if $\sum_{g \in I_i} P_g < \sum_{d \in I_i} P_d$, then do $D_{I_i}^s < \sum_{g \in I_i} P_g$ in stage $s$.
  - if $\sum_{g \in I_i} P_g > \sum_{d \in I_i} P_d$, then do $D_{I_i}^s < \sum_{d \in I_i} P_d$ in stage $s$.
  
  b) for each subnet $I_i$ without generators $g \in I_i$; do $D_i^s = 0$.

Here, $P_g$ is the power of the generator $g$, $P_d$ is the power demand and $D_{I_i}^s$ is the demand on the island $I_i$ during stage $s$.

**Step 8.** *Satisfied demand*: calculate (12),

$$SD = \frac{\sum_{i \in I} D_{I_i}^s}{SD_{base}} \; in \; iteration \; s; \; Unit : [p.u] \qquad (12)$$

**Step 9.** *Iterations*: iterate $s = s+1$ and go to Step 4.
**Step 10.** *End*: if $|P_k^s| < P_k^{max} \; \forall k$, the algorithm ends.

The first step in this process is to gather technical information about the electrical network and calculate power flows to ascertain the lines' maximum transfer capacity. Then, after determining the changes in the flows and removing any overloaded power lines brought on by the redistributed network flows, it randomly eliminates one asset. Keep in mind that unpredictable events, including unintentional human mistake or technical issues with gear and equipment, might start cascading failures. The propagation of the cascading failure may cause the protective mechanisms in the power lines to trip repeatedly, which can result in the construction of various islands in the system. To identify the subsets created during the disintegration phases, this technique thus uses a transversal graph algorithm. Here, the issue is solved more simply by using the Depth-First Search (DFS) method [43]. As a result, the DFS algorithm recognizes and arranges each island for a proper simulation each time one or more power lines are cut during the disintegration stage. Similar to other islands, these ones must maintain a balance between generation and demand, hence load shedding is used to do so. During the disintegration process, isolated components or subnets without creation are regarded as unsatisfied loads. The iterative process continues until all assets are segregated or there are no overloaded parts.

*2.3. Resilience procedure*
Resilience is dependent on both the robustness and the speedy recovery of the disconnected load because, after a significant disruption, robustness determines how quickly an electrical system degrades. Due to this, a mixed-integer optimization problem is suggested to restore the system's connection and loads after a cascading failure. The quantified ideal resilience characteristic and the transmission line states are the outputs of the optimization process. The system's resilience is represented by the RD index for demonstration purposes. Maximizing this resilience measure after the cascading failure depicted by Algorithm 2 is the optimization goal. Several limitations apply to this optimization issue. The iterative process used to identify the power lines that must be shut down at each step of a power system recovery is described in Algorithm 3.

In essence, Algorithm 3 initializes both the recovered demand and the initial operational condition of the power lines using the ultimate disintegration of a power system as input data. The maximum number of lines that may be reconnected as well as the deployment of generators at each stage of recovery are also taken into account. After that, it develops an optimization problem based on the DC power flow equations and determines the lowest and maximum thresholds for each equation. Algorithm 2 determines the power lines' maximum threshold. The goal function for the relevant recovery iteration is maximized once the set of equations has been formed. The restored demand and the power lines that need to be shut down during the restoration iteration make up the output. The new closed power lines are then put in their appropriate equation once these findings have been saved. The algorithm terminates if all power lines are closed; otherwise, it continues the process until all remaining lines are closed.

**Algorithm 3** Resilience
**Input:** outputs of Algorithm 2, i.e., Satisfied Demand (*SD*) in the last disintegration stage (*s*), set of sub-nets (*I*), set of isolated elements (*E*) and state of the links ($\mu_k$). Likewise, the maximum number of lines to reconnect ($N_c$) in each iteration (*r*).
**Output:** recovery of the electrical power system. Recovered demand (*RD*) and state of the links ($\mu_k$) in each iteration (*r*).
**Step 1.** *Start*: initialize $RD_r = SD_{sfinal}$ and $\mu_k=1$ for closed lines and $\mu_k=0$ for open lines at $r=1$. The initial satisfied demand and the states of the lines correspond to the final disintegration state obtained with Algorithm 2.
**Step 2.** *Optimization problem*: consider (13) to (21),

$$\max \left( RD_r - RD_{r-1} \right) \qquad (13)$$

Subject to:

$$P_g^{\min} \leq P_g^r \leq P_g^{\max}, \quad \forall g \in G \qquad (14)$$

$$P_k^{\min} \cdot \mu_k^r \leq P_k^r \leq P_k^{\max} \cdot \mu_k^r, \quad \forall k \in K \qquad (15)$$

$$\Delta_n^{\min} \leq \Delta_n^r \leq \Delta_n^{\max}, \quad \forall n \qquad (16)$$

$$-\sum_{k \in K} P_k^r - \sum_{g \in G} P_g^r - \sum_{d \in D} P_d^r = 0, \quad \forall n \qquad (17)$$

$$B_k \cdot \left( \Delta_n^r - \Delta_m^r \right) - P_k^r \geq 0, \quad \forall k \qquad (18)$$

$$-B_k \cdot \left( \Delta_n^r - \Delta_m^r \right) - P_k^r \leq 0, \quad \forall k \qquad (19)$$

$$\sum_{k \in K} \mu_k^r \leq N_c \qquad (20)$$

$$RD_r = \sum P_n^r, \quad \forall n \qquad (21)$$

The maximum thresholds in (15) are initially determined in (11).
**Step 3.** *Solve the optimization problem*: maximize (13), subject to constraints (14)-(21) in each recovery step $r$.
**Step 4.** *Solution*: save the results of $RD_r$ and $\mu_k^r$; set the variables $\mu_k^r$ restored as constants $\mu_k^r=1$ for all subsequent iterations.
**Step 5.** *Evaluation*: if $\forall k\in(K-k')$: $\mu_k^r=1$ go to Step 7; otherwise, go to Step 6.
**Step 6.** *Iterations*: iterate $r=r+1$ and go to Step 3.
**Step 7.** *End*: if $\forall k\in(K-k')$: $\mu_k^r=1$; the algorithm ends.
**Step 8.** *Energy not supplied*: calculate the ENS index for the resilience curve, i.e. the area above the curve.

### III. Case Studies

This section explains the IEEE 24-bus RTS, which is used to build eight case studies with various topologies [38]. In other words, the original system is utilized while more lines are added to create other networks for comparison. Case studies are offered first, followed by a description of the procedures followed throughout the reliability, robustness, and resilience simulations.

*3.1. Test system*
The improved two-area IEEE 24-bus RTS with three wind power generators is shown in Fig. 2. (WPP). This network consists of 38 power lines and transformers, 24 buses, and 33 generators. 2850 MW is its maximum demand. You may refer to another source [38] for information on the parameters of the lines, the load characteristics, and the input data for the stochastic failure model for buses, transformers, and lines. The scientific literature has several references to this large-scale test network.

Assessments of dependability, robustness, and resilience are suggested for eight distinct case studies based on the prior system. The case studies include taking the existing system and merging and adding three more power lines to it (14-15, 14-24, and 6-9). The goal is to have many topologies of the same system so that comparisons may be made between them. According to the guidelines given in [44], the lines added to the original system fulfill type n-1 and n-2 contingencies. In light of the test system depiction in Fig. 2, the eight case studies are as follows:

> **Case 1:** the original system.
> **Case 2:** addition of line 14-15 to the original power system.
> **Case 3:** addition of line 14-24 to the original power system.
> **Case 4:** addition of line 6-9 to the original power system.
> **Case 5:** addition of lines 14-15 and 14-24 to the original power system.
> **Case 6:** addition of lines 14-15 and 6-9 to the original power system.
> **Case 7:** addition of lines 14-24 and 6-9 to the original power system.
> **Case 8:** addition of lines 14-15, 14-24 and 6-9 to the original power system.

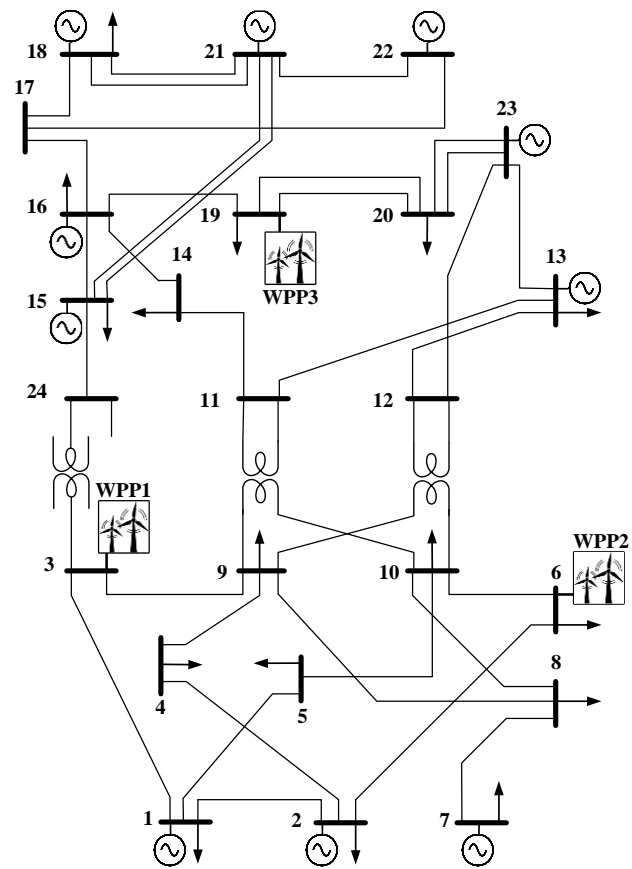

**Fig. 2. Topology of the IEEE 24-bus RTS.**

*3.2. Simulation guidelines in the analysis*

In order to provide a thorough and accurate examination of the many variables under study, different criteria were followed while using Algorithms 1, 2, and 3 in the eight case studies mentioned above. On the one hand, dependability assessment is a traditional power system analysis, therefore our study followed the general framework of previously published studies in this area of study. In each network, 1500 one-year iterations were run, yielding covariance values that were lower than 6% in every instance [45].

In contrast, robustness assessment is a more involved process that takes into account several system factors and features. For instance, the resilience of an electrical system might vary based on the location of the first breakdown, the degree of line congestion, the operating assets, and the load level, among other factors. Because of this, some studies choose to define resilience in terms of topology and structure [46] in order to make it independent of other network-related variables. Naturally, the benefit of this kind of study is that it provides a different view of the system. In this study, it was proposed to start the network disintegration process by removing the lines next to the buses (except from buses 6, 9, 14, 15, and 24 since new lines were added) in order to conduct a thorough assessment of resilience.

Furthermore, given that the system had a constant load, it was also suggested to consider different levels of overload in the links for the same scenario to obtain different states of disintegration for the initiating event. As a result, a total of 114 scenarios were executed with α = 1, 1.1, 1.2, 1.3, 1.4 and 1.5 in each case study, that is, a total of 912 simulations for the eight cases. Finally, the robustness of each case was measured by averaging the set of results obtained, which provided an overview of system performance.

On the other hand, the robustness technique was used before doing the resilience assessment using the eight examples' mean states of disintegration. In each recovery stage, a maximum of three lines were thought to be the most that could change states. The number of lines that may be used to repair an electrical system that has failed relies on the network's physical features and the protocols used by each control center. Three electricity lines were solely used in this paper for simulation. In the end, the ENS index was computed for each example, assuming that the time required to plan, carry out, and check each period of re-dispatch and reconfiguration took on average around 15 minutes.

### IV. Simulation Results

This section discusses the simulation results obtained after performing the reliability, robustness, and resilience evaluations in the eight case studies described above. The three procedures were programmed in MATLAB R2021a and were executed in a personal computer with a 3.40 GHz Intel® Core TM i7 CPU and 16 GB of RAM. The run times for the reliability, robustness and resilience studies were 294.99 h, 167.91 s and 31.81 s, respectively.

First, Table 1 shows the different reliability indicators calculated for the eight cases after applying Algorithm 1, and Fig. 3 presents the convergence results of the EENS indicator. The best-case corresponded to Case 6, where an improvement of 2.19% was obtained compared to Case 1. However, it was also observed that Cases 2 and 7 had EENS values very close to each other, even though Case 7 was more connected than the other. The latter was because line 14-15 reduced the congestion of the two lines adjacent to bus 14. In fact, this line also coincided with Case 6. In Case 8, the improvement was 1.86% over the original system, indicating an improvement in the system performance. The rest of the indicators had similar behaviors to the analysis carried out previously. The results showed that reliability was improved by adding lines to the original system; however, certain lines located on buses with poor connectivity showed better results. From highest to lowest reliability are Cases 6, 5, 8, 2, 7, 4, 3 and 1.

Second, Fig. 4 depicts the dispersion of the last SD value after applying Algorithm 2, as well as the mean values obtained in the eight cases studied. The mean value of the SD index was 0.34 p.u, 0.38 p.u, 0.34 p.u, 0.35 p.u, 0.40 p.u, 0.39 p.u, 0.35 p.u and 0.40 p.u for each case, respectively. The plotted results show that the cases had different satisfied demand values at the end of network collapse, indicating that the redistributed flows after the initial disturbance were different in each case. However, when averaging the sets of results for each case, the robustness of Case 8 improved by 9.43% compared to Case 1. That is, the most connected case was the most robust to cascading failures. The latter makes sense since the power lines of the system were less congested. The results also showed that all the cases where one or two lines were added improved the robustness of the original system. Cases 3 and 4 had improvements of 0.52% and 1.17% when considering less-connected cases compared to Case 1, respectively. However, Case 2 improved by 5.71% compared to the original system and was even better than Case 7 with two lines. Note that Case 2 corresponds to the addition of line 14-15, which was also identified as an asset that improves the system's reliability. The robustness of cases can be ordered from highest to lowest as Cases 8, 5, 6, 2, 7, 4, 3 and 1. These findings suggest that, in general, more meshed topologies are more robust against cascading failures than less meshed topologies.

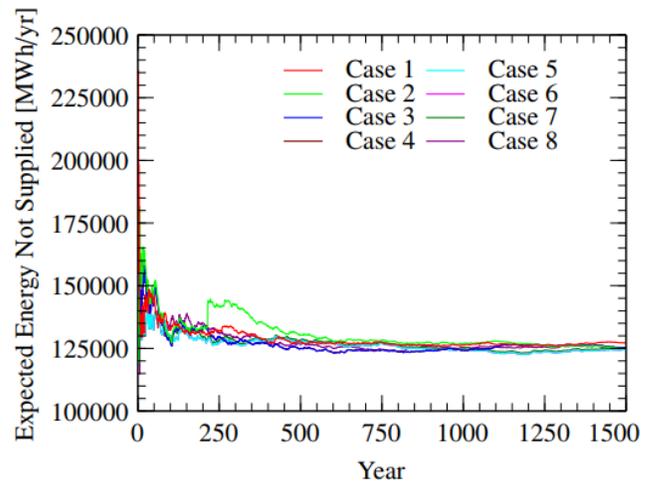

**Fig. 3. Convergence of the EENS indicator for the cases.**

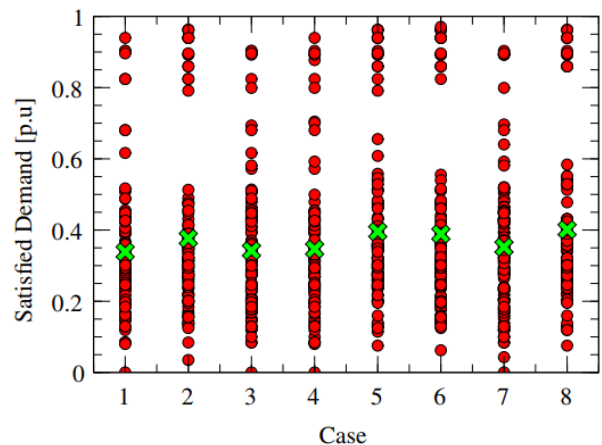

**Fig. 4. Dispersion of the robustness results for the eight cases studied. The red circles represent the final SD values in each simulation, and the green crosses represent the mean values obtained in each case.**

Third, based on the average resilience of the eight case studies, Fig. 5 displays the best recovery curves. For enhanced realism, the scenarios began with a variety of fulfilled demand values as well as disintegration states. Cases 1 and 5 originated from topologies with 27 and 32 missing power lines, respectively, to put a number on this.

Table 1. Reliability results

| Index | Case 1 | Case 2 | Case 3 | Case 4 | Case 5 | Case 6 | Case 7 | Case 8 |
|---|---|---|---|---|---|---|---|---|
| EENS (MWh/yr) | 126985.4 | 124874 | 125351 | 125349 | 124522 | 124508 | 125029 | 124921 |
| EDNS (MW) | 15.23 | 14.56 | 14.68 | 14.68 | 14.22 | 14.22 | 14.31 | 14.25 |
| EFLC (outage/yr) | 18.99 | 18.71 | 18.92 | 18.91 | 18.99 | 18.98 | 18.80 | 19.02 |
| LOLE (h/yr) | 732.68 | 719.36 | 730.25 | 730.26 | 720.18 | 720.21 | 721.05 | 724.16 |
| LOLP (%) | 7.94 | 7.79 | 7.85 | 7.84 | 8.22 | 8.24 | 8.24 | 8.31 |
| ADLC (h/outage) | 39.32 | 39.04 | 39.01 | 39.02 | 38.05 | 38.04 | 38.25 | 38.14 |

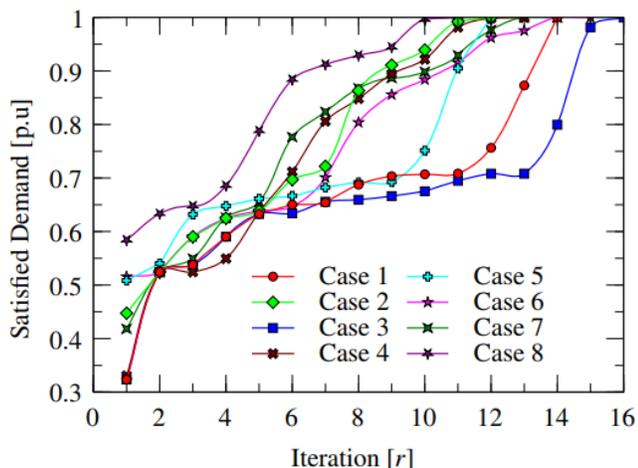

Fig. 5. Optimal recovery curves for the eight cases studied. The plotted curves represent the total satisfied demand within the system. These values are obtained by adding the RD index after each set of maneuvers.

The curves in Fig. 5 illustratively test the general concept and demonstrate the advantage of the restoration strategy proposed in Algorithm 3. The results showed that each of the systems recovered all the disrupted loads; however, some of them were superior because they required fewer maneuvers to restore the load more quickly. The ENS index was 3277.82 MWh, 2108.23 MWh, 3701.65 MWh, 2233.49 MWh, 2564.52 MWh, 2371.71 MWh, 2110.90 MWh and 1359 MWh for each case, respectively. Considering this last indicator, Case 8 was the most resilient because it improved by 58.53% compared to Case 1, while Case 3 was the least resilient because it worsened by 12.93% compared to Case. 1. Cases 2 and 7 improved by 35.68% and 35.60% compared to the original system, which placed them in the second and third positions. The cases can be ordered from highest to lowest resilience as Cases 8, 2, 7, 4, 6, 5, 1 and 3. The resilience evaluation showed that the network topology influenced the recovery of the system. For example, Cases 1 and 3 had different recovery curves even though they started from very similar values of satisfied demand. In other words, topology plays a fundamental role in the design of resilient systems.

Table 2 provides a more detailed look at the findings by displaying the EENS, SD, and ENS indices' improvement percentages for the reliability, robustness, and resilience assessments. In comparison to the original system, these figures are given as percentage increases or decreases (Case 1). A three-dimensional depiction of these findings is shown in Fig. 6.

On the one hand, the findings demonstrated that, with the exception of Case 3, which had the lowest performance in terms of resilience, the majority of the topologies improved on dependability, robustness, and resilience. Similarly, Case 8's topology was the most durable and resilient of all the systems, although it was somewhat less dependable than Cases 5 and 6's topologies due to the additional 0.33% of energy that was not delivered. In terms of the three qualities, Cases 2, 5, and 6's topologies were intermediate. The topologies of Cases 4 and 7 performed well in terms of reliability and resilience, but they performed poorly in terms of robustness since their improvement percentages were small in comparison to the topologies that were more robust. The performance of these two topologies, however, was better than the Case 3 topology.

Table 2. Percentages of increase of the SD, EENS, and ENS indicators in the robustness, reliability, and resilience evaluations compared to Case 1

| Cases | Robustness [$\Delta$SD (%)] | Reliability [$\Delta$EENS (%)] | Resilience [$\Delta$ENS (%)] |
|---|---|---|---|
| Case 1 | 0.00 | 0.00 | 0.00 |
| Case 2 | 1.825 | 5.732 | 35.662 |
| Case 3 | 1.593 | 0.623 | 13.93 |
| Case 4 | 1.594 | 1.184 | 31.841 |
| Case 5 | 2.215 | 8.623 | 21.742 |
| Case 6 | 2.231 | 7.581 | 27.864 |
| Case 7 | 1.778 | 2.038 | 35.625 |
| Case 8 | 1.864 | 9.473 | 58.691 |

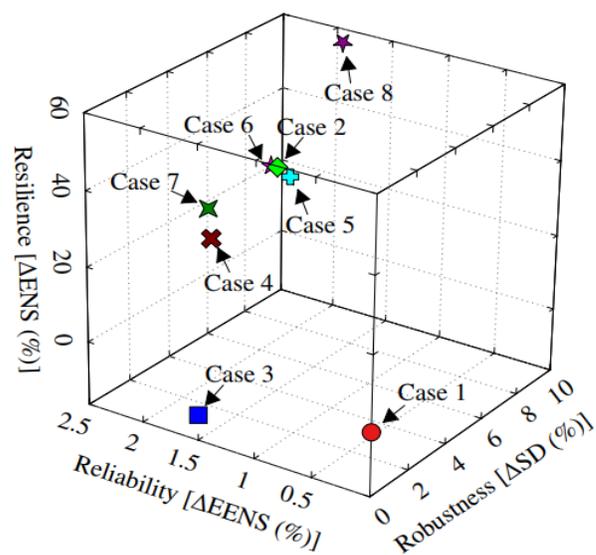

Fig. 6. Results integrated into an R3 concept. The plotted values represent the increases or decreases in the EENS, SD and ENS indices compared to Case 1.

Adding lines allowed for an increase in the power transfer capacity between various zones and decreased congestion in the power lines, which in turn had a favorable influence on the operating conditions of the original system, according to the findings included into an R3 concept. Additionally, they aided in adjusting to various network disruptions that were simulated and permitted the best resource management throughout the recovery phase. To make the R3 idea a better compromise for the construction of electrical power systems, it is crucial to remember that certain lines were more useful than others.

## V. Conclusion

In order to study a power system's dependability, robustness, and resilience from an integrated point of view, this work provided a methodological framework based on data. Reliability was assessed using a sequential Monte Carlo method, robustness was measured using a cascading failure mechanism, and resilience was calculated using a recovery method based on a mixed-integer optimization problem. Eight case studies based on the well-known IEEE 24-bus RTS were utilized for this investigation, and several measures of reliability, robustness, and resilience were determined. The findings were extensively described in a three-dimensional representation that took into account the ranking of each example in each idea and were shown both visually and quantitatively.

The benefit of integrating the three ideas rather than showing them separately was shown by this illustration. The results demonstrated that although it is not always feasible to ensure that a power system's most mesh topology is the best from the perspective of each criteria, in general, it provides the greatest outcomes for supply security. These findings may be used by a power system planner to decide which power grid topology investments are the best.

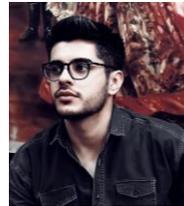

**Farhad Samadi Gazijahani** (M'20) received the PhD Degree in Electric Power Engineering from Azarbaijan Shahid Madani University, Tabriz, Iran, in 2021, graduating with first class honors. He was selected as a 2% of scientists in the world. Moreover, he is a member of *Iran's National Elites Foundation*. His main research interests include power system operation and planning, power system resilience, electricity markets, and optimization techniques.

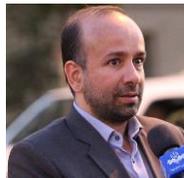

**Rasoul Esmaeilzadeh** (M'18) received the Ph.D. degree from the Azarbaijan Shahid Madani University, Tabriz, Iran. He joined the Azarbaijan Reginal Electric Company. His main research interests include power system reliability and security, renewable energy, power system economics and energy market.